\documentclass[a4]{article}

\usepackage{amsfonts}
\usepackage{bm}
\usepackage{bbm}
\usepackage{graphicx}

\newcommand{\be}{\boldsymbol{e}}
\newcommand{\bh}{\boldsymbol{h}}
\newcommand{\bl}{\boldsymbol{l}}

\newcommand{\bu}{\boldsymbol{u}}
\newcommand{\bx}{\boldsymbol{x}}
\newcommand{\by}{\boldsymbol{y}}

\newcommand{\bN}{\boldsymbol{N}}

\newcommand{\bX}{\boldsymbol{X}}

\newcommand{\bpi}{\boldsymbol{\pi}}
\newcommand{\btheta}{\boldsymbol{\theta}}
\newcommand{\bxi}{\boldsymbol{\xi}}
\newcommand{\R}{{\mathbb R}}
\newcommand{\Z}{{\mathbb Z}}

\begin{document}
\title{Quasistationarity and extinction for population processes}
\author{Damian Clancy\\
Department of Actuarial Mathematics and Statistics\\ Maxwell Institute for Mathematical Sciences\\ Heriot-Watt University\\ Edinburgh\\ EH14 4AS\\ UK\\
d.clancy@hw.ac.uk}
\maketitle

\begin{abstract}
We consider stochastic population processes that are almost surely absorbed at the origin within finite time.
Our interest is in the quasistationary distribution,~$\bu$, and the expected time, $\tau$, from quasistationarity to extinction, both of which we study via WKB approximation.
This approach involves solving a Hamilton-Jacobi partial differential equation specific to the model.
We provide conditions under which analytical solution of the Hamilton-Jacobi equation is possible, and give the solution.
This provides a first approximation to both~$\bu$ and~$\tau$.
We provide further conditions under which a corresponding `transport equation' may be solved, leading to an improved approximation of~$\bu$.
For multitype birth and death processes, we then consider an alternative approximation for~$\bu$ that is valid close to the origin, provide conditions under which the elements of this alternative approximation may be found explicitly, and hence derive an improved approximation for~$\tau$.
We illustrate our results in a number of applications.
\end{abstract}



\section{Introduction}\label{introduction}
For many biological population models, including models from ecology~\cite{MV12,OM10}, epidemiology~\cite{N91,N99,AB00,N11,C18b} and immunology~\cite{SLBHM10}, eventual extinction (either of the whole population or of some sub-population of interest) is certain, but the population can settle to an apparent equilibrium for a long time before extinction occurs.
Two particular objects of interest are then (i)~the quasistationary  distribution,~$\bu$, that the process settles to prior to eventual extinction; and (ii)~the time taken for extinction to occur.
The time from quasistationarity to extinction is in general an exponentially distributed random variable~\cite{DP13}, so that its distribution is fully determined by its expected value~$\tau$.
Provided the relevant quasistationary distribution is known to exist, it is straightforward to write down equations satisfied by $\bu$ and $\tau$, equations~(\ref{QSD}) and~(\ref{tau_formula}) below, but in general far from straightforward to evaluate their solutions.
Approximation methods are therefore of great interest.

We consider a sequence of Markov processes $\{ \bX^{(N)} (t) : t \ge 0 \}$ on~$\Z_+^k$ indexed by a parameter~$N$ representing `typical' population size, such that the origin is an absorbing state, the remainder of the state space consists of a single communicating class~$C$, and absorption at the origin occurs within finite time with probability~1.
This sequence of processes is assumed to be density dependent in the sense of chapter~11 of~\cite{EK05}, so that under mild technical conditions~(theorem 11.2.1 of~\cite{EK05}), the scaled processes $\bX^{(N)}(t)/N$ may be approximated by the solution of a system of ordinary differential equations, system~(\ref{ODEsystem}) below.
We restrict attention to the situation where this deterministic system has a unique stable equilibrium point $\by^*$ in the interior of $\R_+^k$, and the origin is an unstable equilibrium point.

A natural strategy, and one widely used in applications~\cite{N99,LZM07,AB00,R06}, is to approximate the quasistationary distribution~$\bu = \{ u_{\bx} : \bx \in C \}$ by a multivariate normal distribution.
This gives a reasonably good approximation close to the mode at $N \by^*$, but the approximation becomes less accurate away from the mode, and, in particular, fails in the tail of the distribution close to the origin---see, for example, figure~3 of~\cite{LZM07} and figure~4 of~\cite{C18b}.
The multivariate normal approximation consequently does not, in general, lead to a good approximation for the expected extinction time~$\tau$.

An alternative approach is to approximate~$\tau$ by the expected hitting time at the origin of an approximating diffusion process, evaluated as the solution of a Kolmogorov backward partial differential equation~\cite{WGPAI14,V21}.
This technique can give a reasonably good approximation for finite~$N$, but is known to fail in the large~$N$ limit~\cite{DSS05,CT18}.

The approach that we will take is based around a WKB (Wentzel, Kramers, Brillouin) approximation~\cite{OM10,AM17}.
The WKB approach supposes that, as $N \to \infty$,  writing $\by = \bx / N$, the elements of the quasistationary distribution $\bu$ may be approximated in the form
\begin{eqnarray}
u_{\bx} &=& M_N \exp \left( - N V ( \by ) - V_0 ( \by ) + O (1/N) \right) \mbox{ for } \bx \in C  ,
\label{ansatz}
\end{eqnarray}
for some functions $V ( \by )$, $V_0 ( \by )$ that do not depend upon $N$, and some $M_N$ that does not depend upon~$\by$.
We then seek functions $V ( \by )$, $V_0 ( \by )$ and constant $M_N$ such that the conditions for $\bu$ to be a quasistationary distribution are indeed satisfied, in the limit as $N \to \infty$.
The WKB approach is sometimes referred to as `semi-rigorous'; as such, we will avoid the terminology of `theorem' and `proof', and instead refer to our `results' and `derivations.'

By neglecting the second order term $V_0 ( \by )$ in the formula~(\ref{ansatz}) and imposing the constraint that $V ( \by )$ be of quadratic form, we recover the multivariate normal approximation.
It is thus apparent that the WKB approximation provides considerably greater flexibility than the multivariate normal approximation.
In applications, the WKB approximation is found to perform considerably better than the multivariate normal distribution in approximating the quasistationary distribution---see, for example, figure~4 of~\cite{C18b} and figure~4 of~\cite{BM11}.

For a variety of 1-dimensional systems, it is possible to find $V ( \by )$, $V_0 ( \by )$ and $M_N$ explicitly~\cite{AM10,AM17}.
The approximation~(\ref{ansatz}) fails in the tail of the distribution close to extinction, since $V_0 ( \by )$ diverges at the origin.
However, an alternative approximation to $\bu$ may be found by linearising  transition rates close to the origin, giving an un-normalised approximation valid in the relevant tail of the distribution, and then normalising by matching with the WKB approximation.
This leads to a relationship of the form
\begin{eqnarray}
\tau &\sim& {K \over \sqrt{N}} \exp ( NA ) ,
\label{general_form}
\end{eqnarray}
where we use~$\sim$ to denote that the ratio of the two sides converges to~1 as $N \to \infty$, and $A$,~$K$ are constants that can be found in terms of the parameters of the process.

For models in $k>1$ dimensions, it is usually only possible to find the leading order term $V ( \by )$ numerically as the solution of the Hamilton-Jacobi partial differential equation~(\ref{HJE}) below.
This leads to results of the cruder form
\begin{eqnarray}
\lim_{N \to \infty} \frac{\ln \tau}{N} &=& A ,
\label{log_tau}
\end{eqnarray}
where the leading order constant~$A$ must be evaluated via numerical solution of a system of ordinary differential equations in $2k$ dimensions~\cite{AM17,NBKFB17,BFB16,DMRH94,KM08}.
An exception is provided by~\cite{C18b}, where a relation of
the form~(\ref{general_form}) and explicit formulae for $A$,~$K$ are obtained for a $k$ dimensional susceptible–infectious–susceptible (SIS) infection model in a heterogeneous population, example~1 below.
Solution is based upon the observation that for the model of~\cite{C18b}, the $k$ dimensional Hamilton-Jacobi equation~(\ref{HJE}) for $V ( \by )$ can be separated into $k$~components, each of which can be straightforwardly solved as an ordinary differential equation, and then the `transport equation'~(\ref{transport}) can similarly be solved for $V_0 ( \by )$.

A rather different approach was taken in~\cite{BC23} to obtain a relationship of the form~(\ref{general_form}) together with explicit formulae for~$A,~K$ for a class of multitype birth and death processes, example~3 below.
A modified version of the process $X^{(N)} (t)$ is considered, re-started in a particular way each time the process hits the origin.
The stationary distribution of this restarted process is found, and used to study the time to extinction of the original process.

In the current paper, we investigate the extent to which the approach of~\cite{C18b} can be extended to a much more general class of models than those considered by~\cite{BC23}.
Our approach relies upon asymptotic (weaker) versions of Kolmogorov's criterion for reversibility~\cite{K11}.
Kolmogorov's criterion states that for an irreducible, positive recurrent Markov process with state space~$S$, if the transition rates $\{ q_{\bx,\bx^\prime} : \bx , \bx^\prime \in S \}$ satisfy
\begin{eqnarray}
q_{\bx_1 , \bx_2} q_{\bx_2 , \bx_3} \cdots q_{\bx_{n-1} , \bx_n} q_{\bx_n , \bx_1}
&=&
q_{\bx_n , \bx_{n-1}} q_{\bx_{n-1} , \bx_{n-2}} \cdots q_{\bx_2 , \bx_1} q_{\bx_1 , \bx_n}
\label{Kolmogorov}
\end{eqnarray}
for every finite sequence of states $\bx_1 , \bx_2 , \ldots , \bx_n \in S$, then the stationary distribution~$\bpi$ of the process satisfies the detailed balance equations
\begin{eqnarray}
\bpi_{\bx} q_{\bx, \bx^\prime} &=& \bpi_{\bx^\prime} q_{\bx^\prime,\bx} \mbox{ for } \bx , \bx^\prime \in S .
\label{detailed_balance}
\end{eqnarray}
The processes $\bX^{(N)} (t)$ that we consider are not irreducible, having an absorbing state at the origin, and a unique stationary distribution that assigns probability~1 to the state $\bx = {\bf 0}$.
However, we shall show that under various weaker versions of the Kolmogorov criterion~(\ref{Kolmogorov}), asymptotic versions of detailed balance hold, allowing us to find explicitly the functions $V ( \by )$, $V_0 ( \by )$ and constant $M_N$ in the asymptotic approximation~(\ref{ansatz}).

In particular, if the process $\bX^{(N)} (t)$ satisfies the Kolmogorov condition~(\ref{Kolmogorov}) for all~$N$, then our conditions~(\ref{Kolmogorov0},\ref{irrotational},\ref{Kolmogorov1},\ref{irrotational_2}) below are all satisfied.
This may be seen by substituting the density dependent form of the transition rates from equation~(\ref{rates}) into the relationship~(\ref{Kolmogorov}) and expanding in powers of $N$.
We find that leading order terms correspond to our condition~(\ref{Kolmogorov0}).
Second order terms around parallelograms (small loops of the form $\bx \to \bx + \bl_1 \to \bx + \bl_1 + \bl_2 \to \bx + \bl_2 \to \bx$) correspond to our condition~(\ref{irrotational}), or condition~(\ref{irrotational_for_BD_process}) in the case of multitype birth and death processes.
Second order terms around more general closed paths imply our condition~(\ref{Kolmogorov1}).
Third order terms around parallelograms correspond to our condition~(\ref{irrotational_2}), or condition~(\ref{irrotational_for_BD_process2}) for multitype birth and death processes.

On the other hand, if we substitute the density dependent form of the transition rates from equation~(\ref{rates}), together with the assumed WKB form~(\ref{ansatz}) for the quasistationary distribution~$\bu$, into the detailed balance equations~(\ref{detailed_balance}) and expand in powers of $N$, we find that leading order terms correspond to our equation~(\ref{linear_system}) with $\btheta = \frac{\partial V}{\partial \by}$, which allows explicit solution of the Hamilton-Jacobi equation~(\ref{HJE}) for $V ( \by )$.
Second order terms then correspond to our equation~(\ref{linear_system2}) with $\btheta^0 = \frac{\partial V_0}{\partial \by}$, which allows explicit solution of the transport equation~(\ref{transport}) for~$V_0 ( \by )$.

Following the statement of our results in section~\ref{results}, in section~\ref{examples} we exhibit examples of biological population models satisfying the relevant asymptotic Kolmogorov criteria
~(\ref{irrotational_for_BD_process},\ref{irrotational_for_BD_process2},\ref{Kolmogorov0},\ref{irrotational},\ref{Kolmogorov1},\ref{irrotational_2}), and for which we can therefore provide explicit results.
There are, of course, many other population processes that do not satisfy these criteria.
It is of particular interest to understand the conditions under which exact solution of equation~(\ref{HJE}) is possible because numerical solution of the Hamilton–Jacobi equation~(\ref{HJE}) otherwise requires the solution of a high-dimensional system of ordinary differential equations (the characteristic equations) subject to boundary conditions at times $t = - \infty$ and $t = + \infty$.
The relevant characteristic curve is maximally sensitive to perturbations, so that numerical solution is usually only feasible in low dimensions~\cite{AM17,BFB16,BM11,NBKFB17}.
Numerical solution techniques generally require an initial guess of the solution curve, and only work well provided the initial guess is close to the true solution.
When a particular population process of interest does not satisfy our conditions~(\ref{Kolmogorov0},~\ref{irrotational}), but is closely related to a process that does, then one option is to use the explicit solution $V ( \by )$ for the process that does satisfy conditions~(\ref{Kolmogorov0},~\ref{irrotational}) as our initial guess when numerically solving to find the characteristic curve for the process of interest.
An example of this approach is given in~\cite{CS23b}.

Our formula~(\ref{tau_BD}) giving the asymptotic form of the time to extinction from quasistationarity is analogous to the Eyring-Kramers formula~\cite{B13} characterising the expected transition time between locally stable equilibrium points for a reversible diffusion process.
For processes on a discrete state space, the recent paper~\cite{JJL21} studies the exit time and exit distribution from a basin of attraction of a locally stable equilibrium point for chemical reaction network models, so that transition rate functions are restricted to be of the form shown on p889 of~\cite{JJL21}.
Theorem~4 of~\cite{JJL21} parallels, in this somewhat different context, our result~4 below.

The remainder of the paper is structured as follows.
In section~\ref{model} we set out our model assumptions.
Statements of our results are given in section~\ref{results}, and some particular applications in section~\ref{examples}.
Derivations of our results are given in section~\ref{derivations}.

{\bf Notation}: for any set~$E \subset \R^k$, we denote by $E^{\circ}$, $\partial E$ the interior and boundary, respectively, of the set~$E$.
We denote by $d( \cdot , \cdot )$ the Euclidean distance between two points in $\R^k$ or between a point and a set in $\R^k$, so that in particular $d ( \by , \partial E )$ represents the distance from the point $\by$ to the boundary of the set $E$.
For $i=1,2,\ldots,k$ we denote by $\be_i$ the unit vector with $i$th element equal to~1.

\section{Model assumptions}\label{model}
Consider a sequence of Markov processes $\left\{ \bX^{(N)} (t) : t \ge 0 \right\}$ on $\Z_+^k$ indexed by a parameter~$N$.
We adopt the following assumptions throughout.
\begin{enumerate}
\item The state space $S$ is either the whole of $\Z_+^k$ for all $N$, or the finite set 
\linebreak
$\left\{ \bx = \left( x_ 1 , x_2 , \ldots , x_k \right) : 0 \le x_i \le N_i \mbox{ for } i=1,2,\ldots,k \right\}$ for some $\bN = \left( N_1 , N_2 , \ldots , N_k \right)$ with $N_1 + N_2 + \cdots + N_k = N$.
In the case of finite state space, writing $f_i = N_i / N$ for $i = 1,2, \ldots ,k$, we assume for simplicity that $f_1 , f_2 , \ldots , f_k$  do not vary with $N$ (so when we write, for instance, $N \to \infty$, it is implicit that $N$ increases through a subsequence of integers such that $N f_i$ is always integer-valued for all~$i$), and that $f_1 , f_2 , \ldots , f_k > 0$. 
\label{state_space}
\item The sequence of processes is density dependent in the sense of chapter~11 of~\cite{EK05}.
That is, writing $\tilde S = \R_+^k$ in the case that $S = \Z_+^k$ and $\tilde S = [ 0 , f_1 ] \times [ 0 , f_2 ] \times \cdots \times [ 0 , f_k ]$ in the case of finite state space, transition rates are of the form
\begin{eqnarray}
P \left( \bX^{(N)} ( t + \delta t ) = \bx + \bl \mid \bX^{(N)} (t) = \bx \right) 
&=& N \beta_{\bl} \left( {\bx \over N} \right) + o ( \delta t )
\nonumber 
\\ &&
\hspace*{10mm} \mbox{ for } \bx \in S,\ \bl \in {\cal L} , \label{rates}
\end{eqnarray} 
for some functions $\beta_{\bl} : \tilde S \to \R_+$, where ${\cal L}$ is a finite set, independent of~$N$, consisting of the possible jumps from each state $\bx \in S$.
\item The set~${\cal L}$ spans $\R^k$, and for every $\bl \in {\cal L}$ we also have $- \bl \in {\cal L}$.
\label{minus_l}
\item For every $\bl \in {\cal L}$ we have $\beta_{\bl} ( {\bf 0} ) = 0$.
\label{absorbing}
\item For all $\by \in \tilde S \setminus \{ {\bf 0} \}$, $\bl \in {\cal L}$, we have $\beta_{\bl} ( \by ) > 0$ unless $y_i = 0$ and $l_i < 0$ for some $i$ or (in the case of finite state space) $y_i = f_i$ and $l_i > 0$ for some~$i$, in which case $\beta_{\bl} ( \by ) = 0$. 
\label{communication}
\item The functions $\beta_{\bl} ( \cdot )$ are continuous on $\tilde S$ and twice differentiable on $\tilde S^\circ$.
\label{twice_diff}
\item The ordinary differential equation system
\begin{eqnarray}
{d \by \over dt} &=& \sum_{\bl \in {\cal L}} \bl \beta_{\bl} ( \by )  \label{ODEsystem}
\end{eqnarray}
has two equilibrium points in $\tilde S$: a locally unstable equilibrium point at $\by = {\bf 0}$ and a locally stable equilibrium point  $\by^* \in \tilde S^\circ$, with all eigenvalues of the Jacobian of~(\ref{ODEsystem}) at $\by^*$ having strictly negative real part.
\label{equilm_assumptions}
\item For each $N$ there exists a unique (proper) limiting conditional distribution $\bu = \{ u_{\bx} : \bx \in C \}$ such that for every initial state $\bx_0 \in C$,
\begin{eqnarray*}
u_{\bx} &=& \lim_{t \to \infty} \Pr \left( \bX^{(N)}(t) = \bx \left| \bX^{(N)} ( 0 ) = \bx_0 , \, \bX^{(N)} (t) \in C \right. \right) \mbox{ for } \bx \in C .
\end{eqnarray*} 
\label{QSD_exists}
\end{enumerate}

Assumption~\ref{absorbing} ensures that the state $\bx = {\bf 0}$ is absorbing, while assumption~\ref{communication} ensures that $C = S \setminus \{ {\bf 0} \}$ forms a single communicating class.
Under mild conditions, by Theorem~11.2.1 of~\cite{EK05}, the scaled processes $\bX^{(N)} (t) / N$ converge almost surely over finite time intervals, as $N \to \infty$, to the solution $\by (t)$ of the system~(\ref{ODEsystem}).
In the case of a finite state space, assumption~\ref{QSD_exists} is automatically satisfied~\cite{DS67}.
In the case of an infinite state space, checking assumption~\ref{QSD_exists} can be far from trivial; we refer the reader to~\cite{DP13,CCM19,CV23} for relevant results and discussion.
By theorem~13 and corollaries~7 and~9 of~\cite{DP13}, it follows from assumption~\ref{QSD_exists} that~$\bu$ is a minimal quasistationary distribution for the process $\bX^{(N)} (t)$, that $\bX^{(N)} (t)$ is almost surely absorbed at ${\bf 0}$ within finite time, and that for the process initiated from quasistationarity, so $\bX^{(N)} (0) \sim \bu$, the time until absorption at the origin is exponentially distributed.

\section{Results}\label{results}
Recall that $\bu$ denotes the quasistationary distribution of the process, and $\tau$ the expected time to extinction given that the state of the process at $t=0$ is distributed according to~$\bu$.
Since our results simplify significantly in the case of multitype birth and death processes, we present this case first, before moving on in section~\ref{general_results} to more general population processes.

\subsection{Multitype birth and death processes}\label{BD_results}

In the case ${\cal L} = \{ \be_i , - \be_i : i=1,2,\ldots,k \}$, we have the following results.

{\bf Result 1.}
For $i , j \in \{ 1,2,\ldots,k \}$, denote by $b_{ij}$, $d_i$ the constants
\begin{eqnarray}
b_{ij} = \left. \frac{\partial \beta_{\be_i}}{\partial y_j} \right|_{\by = {\bf 0}}  , &&
d_i = \left. \frac{\partial \beta_{-\be_i}}{\partial y_i} \right|_{\by = {\bf 0}} , \label{bi_di}
\end{eqnarray}
and suppose that
\begin{eqnarray}
d_i > 0 \mbox{ for } i=1,2,\ldots,k,
\mbox{ and for each $i$, there exists some $j$ with } b_{ij} > 0 .
\label{b_d_assumptions}
\end{eqnarray}

Define the vector field $\bh ( \by ) = \left( h_1 ( \by ) , h_2 ( \by ) , \ldots , h_k ( \by ) \right)$ to have components
\begin{eqnarray}
h_i ( \by ) &=& \ln \left( \frac{\beta_{-\be_i} (\by)}{\beta_{\be_i} (\by)} \right) \mbox{ for } i=1,2,\ldots,k. \label{h_components_BD}
\end{eqnarray}

If for all $i,j \in \{ 1,2,\ldots,k \}$ and $\by \in \tilde S^\circ$ we have 
\begin{eqnarray}
{\partial h_i \over \partial y_j} &=& {\partial h_j \over \partial y_i} , \label{irrotational_for_BD_process}
\end{eqnarray}
then 
\begin{eqnarray}
\frac{\ln \tau}{N} &\to&  \int_{\Gamma_0} \bh ( \by^\prime ) \cdot d\by^\prime \mbox{ as } N \to \infty , \label{log_limit}
\end{eqnarray}
where $\Gamma_0$ is any path from $\by^*$, the locally stable equilibrium point of the system~(\ref{ODEsystem}), to ${\bf 0}$ that lies entirely within~$\tilde S$, and the integral is independent of the particular path $\Gamma_0$.

{\bf Result 2.}
With functions $h_i ( \by )$ given by equation~(\ref{h_components_BD}), denote by $\Sigma$ the matrix with entries
\begin{eqnarray}
\sigma_{ij} &=& \left. \frac{\partial h_i}{\partial y_j} \right|_{\by = \by^*} \mbox{ for } i,j = 1,2,\ldots,k . \label{Sigma_matrix_BD}
\end{eqnarray}
Define the vector field $\bh^0 ( \by ) = \left( h^0_1 ( \by ) , h^0_2 ( \by ) \ldots , h^0_k ( \by ) \right)$ to have components
\begin{eqnarray}
h^0_i ( \by ) &=& {1 \over 2} {\partial \over \partial y_i} \ln \left( \beta_{\be_i} ( \by ) \beta_{-\be_i} ( \by ) \right) \mbox{ for } i=1,2,\ldots,k. \label{h0_BD}
\end{eqnarray}

If, in addition to condition~(\ref{irrotational_for_BD_process}), we have that for all $i,j \in \{ 1,2,\ldots,k \}$ and $\by \in \tilde S^\circ$,
\begin{eqnarray}
{\partial h^0_i \over \partial y_j} 
&=& 
{\partial h^0_j \over \partial y_i} ,
\label{irrotational_for_BD_process2}
\end{eqnarray}
then the body of the quasistationary distribution may be approximated as
\begin{eqnarray}
u_{\bx} 
&=& 
\sqrt{\frac{\mbox{det} ( \Sigma )}{(2 \pi N)^k}
\prod_{i=1}^k  
\frac{\beta_{\be_i} \left( y_1 , \ldots , y_{i-1} , y_i^* , \ldots , y_k^* \right) \beta_{-\be_i} \left( y_1 , \ldots , y_{i-1} , y_i^* , \ldots , y_k^* \right)} 
{\beta_{\be_i} \left( y_1 , \ldots , y_i , y_{i+1}^* , \ldots , y_k^* \right) \beta_{-\be_i} \left( y_1 , \ldots , y_i , y_{i+1}^* , \ldots , y_k^* \right)}} 
\nonumber \\ && {} \times
\exp \left( - N  \int_{\Gamma} \bh ( \by^\prime ) \cdot d\by^\prime + O (1/N) \right)
\label{QSD_body}
\end{eqnarray}
for $\bx \in S^\circ$ with $d ( \bx , \partial S ) \ne o(N)$, where $\by = \bx / N$ and $\Gamma$ is any path from $\by^*$ to $\by$ that lies entirely within~$\tilde S$, the integral being independent of the particular path~$\Gamma$.

{\bf Result 3.}
If, in addition to conditions~(\ref{b_d_assumptions}),~(\ref{irrotational_for_BD_process}) and~(\ref{irrotational_for_BD_process2}), we have, recalling the definitions~(\ref{bi_di}), that for all $i,j \in \{ 1,2,\ldots,k \}$,
\begin{eqnarray}
b_{ij}
&=&
b_{ii} , \label{linear_Kolmogorov}
\end{eqnarray}
then denoting by $D$ the unique positive solution of
\begin{eqnarray}
\sum_{i=1}^k \frac{b_{ii}}{D + d_i} &=& 1 ,
\label{D_definition}
\end{eqnarray}
we have
\begin{eqnarray}
\tau &\sim&
\frac{1}{D}
\sqrt{\frac{2 \pi}{N \mbox{det} ( \Sigma )}
{b_{kk}
\prod_{i=1}^k \left. \frac{\partial \beta_{-\be_i}}{\partial y_i} \right|_{\by = (0 , \ldots , 0 , y^*_{i+1} , \ldots , y^*_k )}
\prod_{i=1}^{k-1} \beta_{\be_i} \left( 0 , \ldots , 0 , y_{i+1}^* , \ldots , y_k^* \right) 
\over 
\prod_{i=1}^k
\beta_{\be_i} \left( 0 , \ldots , 0 , y_i^* , \ldots , y_k^* \right) \beta_{-\be_i} \left( 0 , \ldots , 0 , y_i^* , \ldots , y_k^* \right)
}}
\nonumber \\ 
&& {} \times
\exp \left( N \int_{\Gamma_0} \bh ( \by^\prime ) \cdot d\by^\prime \right)
\mbox{ as } N \to \infty , \label{tau_BD}
\end{eqnarray}
where $\Gamma_0$ is any path from $\by^*$ to ${\bf 0}$ that lies entirely within~$\tilde S$, the integral being independent of the particular path~$\Gamma_0$.

{\bf Remark 1.}
It is not immediately apparent that the formula~(\ref{tau_BD}) is invariant under re-labelling of the coordinate axes.
That this is indeed the case under the conditions of result~3 follows from the derivation of the result (section~\ref{derivations}), and becomes clear in specific applications, see section~\ref{examples} below.

\subsection{General population processes}\label{general_results}
{\bf Result 4.}
Suppose that for $a_1 , a_2 , \ldots , a_n \in \Z$ and $\bl_1 , \bl_2 , \ldots , \bl_n \in {\cal L}$ we have
\begin{eqnarray}
\sum_{i=1}^n a_i \bl_i = {\bf 0} 
&\Rightarrow&
\prod_{i=1}^n \left( \frac{\beta_{-\bl_i} ( \by )}{\beta_{\bl_i} ( \by )} \right)^{a_i} = 1 \mbox{ for all } \by \in \tilde S^\circ ,
\label{Kolmogorov0}
\end{eqnarray}
and that
\begin{eqnarray}
\bl_1^T \frac{\partial}{\partial \by} \ln \left( \frac{\beta_{-\bl_2} ( \by )}{\beta_{\bl_2} ( \by )} \right)
&=&
\bl_2^T \frac{\partial}{\partial \by} \ln \left( \frac{\beta_{-\bl_1} ( \by )}{\beta_{\bl_1} ( \by )} \right)
\mbox{ for every } \bl_1 , \bl_2 \in {\cal L} .
\label{irrotational}
\end{eqnarray}

Then the system of equations
\begin{eqnarray}
\bl^T \btheta ( \by ) &=& \ln \left( \frac{\beta_{-\bl} ( \by )}{\beta_{\bl} ( \by )} \right) \mbox{ for all } \bl \in {\cal L}
\label{linear_system}
\end{eqnarray}
has a unique solution $\btheta ( \by )$ for each $\by \in \tilde S^\circ$, and provided the integral converges,
\begin{eqnarray}
\frac{\ln \tau}{N} &\to& \int_{\Gamma_0} \btheta ( \by^\prime ) \cdot d \by^\prime \mbox{ as } N \to \infty ,
\label{log_limit_general}
\end{eqnarray}
where $\Gamma_0$ is any path from $\by^*$ to ${\bf 0}$ that lies entirely within $\tilde S$, and the integral is independent of the particular path~$\Gamma_0$.

{\bf Result 5.}
Suppose that conditions~(\ref{Kolmogorov0}) and~(\ref{irrotational}) hold.
With $\btheta ( \by )$ denoting the solution of equations~(\ref{linear_system}), define the matrix $\Sigma$ to be
\begin{eqnarray}
\Sigma
&=& \left. \frac{\partial \btheta}{\partial \by} \right|_{\by = \by^*} .
\label{Sigma_matrix}
\end{eqnarray}

Suppose that for $a_1 , a_2 , \ldots , a_n \in \Z$ and $\bl_1 , \bl_2 , \ldots , \bl_n \in {\cal L}$ we have
\begin{eqnarray}
\sum_{i=1}^n a_i \bl_i = {\bf 0} 
&\Rightarrow&
\sum_{i=1}^n a_i \bl_i^T \frac{\partial}{\partial \by} \ln \left( \beta_{-\bl_i} ( \by ) \beta_{\bl_i} \right) = 0
\mbox{ for all } \by \in \tilde S^\circ ,
\label{Kolmogorov1}
\end{eqnarray}
and that for every $\bl_1 , \bl_2 \in {\cal L}$,
\begin{eqnarray}
\bl_1^T \left( \frac{\partial^2}{\partial \by^2} \ln \left( \beta_{-\bl_1} ( \by ) \beta_{\bl_1} ( \by ) \right) \right) \bl_2 
&=&
\bl_2^T \left( \frac{\partial^2}{\partial \by^2} \ln \left( \beta_{-\bl_2} ( \by ) \beta_{\bl_2} ( \by ) \right) \right) \bl_1 . 
\label{irrotational_2}
\end{eqnarray}

Then the system of equations
\begin{eqnarray}
\bl^T \btheta^{0} ( \by ) 
&=&
\frac{1}{2} \bl^T \frac{\partial}{\partial \by} \ln \left( \beta_{-\bl} ( \by ) \beta_{\bl} ( \by ) \right)
\mbox{ for all } \bl \in {\cal L}
\label{linear_system2}
\end{eqnarray}
has a unique solution $\btheta^{0} ( \by )$ for each $\by \in \tilde S^\circ$, and the body of the quasistationary distribution may be approximated as follows.
For $\bx \in S^\circ$ with $d ( \bx , \partial S ) \ne o(N)$, writing $\by = \bx / N$, then $\by$ may be expressed (non-uniquely) as $\by = \by^* + \sum_{i=1}^n a_i \bl_i$ for some $a_1 , a_2 , \ldots , a_n \in \R$ and $\bl_1 , \bl_2 , \ldots , \bl_n \in {\cal L}$ such that $\by^* + \sum_{i=1}^j a_i \bl_i \in \tilde S$ for $j=1,2,\ldots,n$, and we then have
\begin{eqnarray}
u_{\bx} &=& 
\sqrt{\frac{\mbox{det} ( \Sigma )}{(2 \pi N)^k}
\prod_{j=1}^n \frac{\beta_{-\bl_j} \left( \by^* + \sum_{i=1}^{j-1} a_i \bl_i \right) \beta_{\bl_j} \left( \by^* + \sum_{i=1}^{j-1} a_i \bl_i \right)}{\beta_{-\bl_j} \left( \by^* + \sum_{i=1}^{j} a_i \bl_i \right) \beta_{\bl_j} \left( \by^* + \sum_{i=1}^{j} a_i \bl_i \right)}}
\nonumber \\ && {} \times
\exp \left( - N \int_{\Gamma} \btheta ( \by^\prime ) \cdot d \by^\prime + O (1/N) \right)
\label{QSD_body_general}
\end{eqnarray}
where $\Gamma$ is any path from $\by^*$ to $\by$ that lies entirely within $\tilde S$, the integral is independent of the particular path~$\Gamma$, and the formula~(\ref{QSD_body_general}) is independent of the particular representation $\by = \by^* + \sum_{i=1}^n a_i \bl_i$ chosen.

\section{Applications}\label{examples}
We defer derivations to section~\ref{derivations}, and first illustrate our results in some specific applications.

{\bf Example 1.}
{\em The susceptible–infectious–susceptible (SIS) infection model of~\cite{C18b} with heterogeneous susceptibilities and infectious periods.} 
This process has finite state space, with $N_i$ representing the (constant) total number of type~$i$ individuals in the population, while
$X_i (t)$ is the number of infectious type~$i$ individuals at time $t$. 
For this model, ${\cal L} = \{ \be_i , - \be_i : i = 1,2,\ldots,k \}$, with transition rate functions $\beta_{\be_i} (\by) = \beta \mu_i \left( f_i - y_i \right) \sum_{j=1}^k y_j$ and $\beta_{-\be_i} (\by) = y_i / \alpha_i$.
Here $\beta > 0$ is an overall infection rate parameter, while for $i=1,2,\ldots,k$, $\mu_i > 0$ gives the level of susceptibility of uninfected type~$i$ individuals and $\alpha_i > 0$ is the mean infectious period of infected type~$i$ individuals. 
Recall also (assumption~\ref{state_space}) that $f_i = N_i / N > 0$ for $i=1,2,\ldots,k$, where $N = N_1 + N_2 + \cdots + N_k$.
Without loss of generality the $\mu_i$ values are scaled so that $\sum_{i=1}^k \mu_i f_i = 1$.
In order that assumption~\ref{equilm_assumptions} is satisfied, we require $\beta \sum_{i=1}^k \alpha_i \mu_i f_i > 1$.

From the definitions~(\ref{bi_di}) we have $b_{ij} = \beta \mu_i f_i$ and $d_i = 1 / \alpha_i$ for $i,j=1,2,\ldots,k$, so that conditions~(\ref{b_d_assumptions}) and~(\ref{linear_Kolmogorov}) are satisfied.
Denoting by $E$ the unique positive solution of
\begin{eqnarray*}
\beta \sum_{i=1}^k \frac{\alpha_i \mu_i f_i}{\alpha_i \mu_i E + 1} &=& 1 ,
\end{eqnarray*}
then from equation~(9) of~\cite{C18b}, the stable equilibrium point $\by^*$ has components
\begin{eqnarray*}
y_i^* &=& \frac{\alpha_i \mu_i f_i E}{1 + \alpha_i \mu_i E} \mbox{ for } i=1,2,\ldots,k.
\end{eqnarray*}
It is straightforward to check that conditions~(\ref{irrotational_for_BD_process}) and~(\ref{irrotational_for_BD_process2}) are satisfied.
Following some algebraic manipulation, including applying the matrix determinant lemma to evaluate $\mbox{det} ( \Sigma )$, the relationship~(\ref{tau_BD}) reduces to
\begin{eqnarray*}
\tau &\sim&
\frac{1}{DE}
\sqrt{\frac{2\pi}{N} \left( \sum_{i=1}^k f_i \left( \frac{\alpha_i \mu_i}{1 + \alpha_i \mu_i E} \right)^2 \right)^{-1}}
\;
\exp \left( N \left(
\sum_{i=1}^k f_i \ln \left( 1 + \alpha_i \mu_i E \right) - \frac{E}{\beta} \right) \right) ,
\end{eqnarray*}
in agreement with formula~(6) of~\cite{C18b} for the case of exponentially distributed infectious periods.

{\bf Example 2.}
{\em The multitype birth and death process with linear birth rates and quadratic death rates described in \cite{CCM19}, Section~2.2.} 
This process has state space $\Z_+^k$ and ${\cal L} = \{ \be_i , - \be_i : i = 1,2,\ldots,k \}$, with transition rate functions $\beta_{\be_i} ( \by ) = \lambda \sum_{j=1}^k y_j$ and $\beta_{-\be_i} ( \by ) = y_i \left( \mu + \kappa \sum_{j=1}^k y_j \right)$ for parameters $\lambda , \mu , \kappa > 0$.
In order that assumption~\ref{equilm_assumptions} is satisfied, we require $k \lambda > \mu$.
Assumption~\ref{QSD_exists} then follows from~\cite{CCM19}.
It is straightforward to check that conditions~(\ref{b_d_assumptions},\ref{irrotational_for_BD_process},\ref{irrotational_for_BD_process2},\ref{linear_Kolmogorov}) are satisfied, and that the relationship~(\ref{tau_BD}) reduces, after some algebraic manipulation, to
\begin{eqnarray*}
\tau &\sim& 
\frac{1}{(k \lambda - \mu )^2} 
\sqrt{\frac{2 \pi \mu \kappa}{N}}
\exp \left( N \left( \frac{k\lambda - \mu}{\kappa} + \frac{\mu}{\kappa}\ln \left( \frac{\mu}{k \lambda} \right) \right) \right) .
\end{eqnarray*}
This is in agreement with formula~(49) and equations~(51) and~(52) of~\cite{BC23}, with the probability~$\omega$ in equation~(52) of~\cite{BC23} being given by $\omega = \mu / k \lambda$, corresponding to the formula immediately before equation~(55) of~\cite{BC23}.

{\bf Example 3.}
{\em The multitype birth and death processes of \cite{BC23}.} 
With ${\cal L} = \{ \be_i , - \be_i : i=1,2,\ldots,k \}$, then setting $\beta_{\be_i} ( \by ) = b_0 \left( \sum_{j=1}^k y_j \right) b_i ( y_i )$ and $\beta_{-\be_i} ( \by ) = d_0 \left( \sum_{j=1}^k y_j \right) d_i ( y_i )$ for appropriate functions $b_0 , b_1 , \ldots , b_k , d_0 , d_1 , \ldots , d_k$ from $\R_+$ to $\R_+$ gives the form of multitype birth and death process considered by~\cite{BC23}.
We require the functions $b_0 , b_1 , \ldots , b_k , d_0 , d_1 , \ldots , d_k$ to be such that our assumptions~\ref{absorbing}--\ref{QSD_exists} are satisfied.
In particular, assumptions~\ref{absorbing} and~\ref{communication} require that for $i=1,2,\ldots,k$, $b_i (0) > 0$ and $d_i (0) = 0$, and that $b_0 (0) = 0$.
In order that condition~(\ref{b_d_assumptions}) is satisfied, we further require that for $i=1,2,\ldots,k$,  $d_i^\prime (0) > 0$, and that $d_0 (0) > 0$ and  $b_0^\prime (0) > 0$.
These conditions correspond to the conditions~(4,5,8) of~\cite{BC23}.
It is then straightforward to check that our conditions~(\ref{irrotational_for_BD_process},\ref{irrotational_for_BD_process2},\ref{linear_Kolmogorov}) are satisfied, and that the relationship~(\ref{tau_BD}) reduces, after some algebraic manipulation, including use of the matrix determinant lemma to evaluate $\mbox{det} ( \Sigma )$, to the result given by formula~(49) together with equations~(41),~(42) and~(55) of~\cite{BC23}

{\bf Example 4.}
{\em Competition processes.}
Taking $k=2$ and ${\cal L} = \{ \be_1 , - \be_1 , \be_2 , - \be_2 , \be_2 - \be_1 , \be_1 - \be_2 \}$ corresponds to the `competition processes' of~\cite{R61}.
Suppose that transition rate functions are of the form
\begin{eqnarray}
\left.
\begin{array}{rcl}
\beta_{-\be_1} ( \by ) &=& a_1 d_0 ( y_1 + y_2 ) d_1 ( y_1 ) c_2 ( y_2 ) d_3 ( y_1 - y_2 ) , \\
\beta_{\be_1} ( \by ) &=& a_2 b_0 ( y_1 + y_2 )  b_1 ( y_1 ) c_2 ( y_2 ) b_3 ( y_1 - y_2 ) , \\
\beta_{-\be_2} ( \by ) &=& a_3 d_0 ( y_1 + y_2 ) c_1 ( y_1 ) d_2 ( y_2 ) b_3 ( y_1 - y_2 ) , \\
\beta_{\be_2} ( \by ) &=& a_4 b_0 ( y_1 + y_2 ) c_1 ( y_1 ) b_2 ( y_2 ) d_3 ( y_1 - y_2 ) , \\
\beta_{\be_1 - \be_2} ( \by ) &=& a_5 c_0 ( y_1 + y_2 ) b_1 ( y_1 ) d_2 ( y_2 ) b_3^2 ( y_1 - y_2 ) , \\
\beta_{\be_2 - \be_1} ( \by ) &=& a_6 c_0 ( y_1 + y_2 ) d_1 ( y_1 ) b_2 ( y_2 ) d_3^2 ( y_1 - y_2 ) ,
\end{array} \right\}
\label{competition_rates}
\end{eqnarray}
for some positive constants $a_1 , a_2 , a_3 , a_4 , a_5 , a_6$, functions $b_0 , b_1 , b_2 , c_0 , c_1 , c_2 , d_0 , d_1 , d_2$ from $\R_+$ to $\R_+$, and functions $b_3 , d_3$ from $\R$ to $\R_+$.
We require these functions to be such that assumptions~\ref{absorbing}--\ref{QSD_exists} are satisfied.
In particular, assumptions~\ref{absorbing} and~\ref{communication} imply that $b_0(0) = d_1 (0) = d_2 (0) = 0$ and that $b_1(0) , b_2(0) , c_1(0) , c_2(0) , b_3(0) , d_3 (0) > 0$.

To check condition~(\ref{Kolmogorov0}), it is more straightforward here to check directly that equations~(\ref{linear_system}) have a unique solution; we find that this is the case provided $a_1 a_4 a_5 = a_2 a_3 a_6$.
It is then straightforward to check that condition~(\ref{irrotational}) is satisfied.
The function $V ( \by )$ given by equation~(\ref{V_general}) may be written as
\begin{eqnarray}
V ( y_1 , y_2 ) &=& 
( y_1^* - y_1 ) \ln \left( \frac{a_2}{a_1} \right) 
+ ( y_2^* - y_2 ) \ln \left( \frac{a_4}{a_3} \right)
+ \int_{y_1 + y_2}^{y_1^* + y_2^*} \ln \left( \frac{b_0(u)}{d_0(u)} \right) du
\nonumber \\ && {} \hspace*{-20mm}
+ \int_{y_1}^{y_1^*} \ln \left( \frac{b_1(u)}{d_1(u)} \right) du
+ \int_{y_2}^{y_2^*} \ln \left( \frac{b_2(u)}{d_2(u)} \right) du
+ \int_{y_1 - y_2}^{y_1^* - y_2^*} \ln \left( \frac{b_3(u)}{d_3(u)} \right) du .
\label{V_competition}
\end{eqnarray}
The relationship~(\ref{log_limit_general}) reduces to $(\ln \tau ) / N \to V ( 0,0 )$ as $N \to \infty$ with $V ( y_1 , y_2 )$ given by equation~(\ref{V_competition}), provided the integrals in equation~(\ref{V_competition}) with $y_1 = y_2 = 0$ all converge.

To check condition~(\ref{Kolmogorov1}), it is more straightforward here to check directly that equations~(\ref{linear_system2}) have a unique solution; we find that this is the case provided that $b_3 (u) d_3 (u) = a_0$ for all $u \in \R$ for some constant $a_0 > 0$.
It is straightforward to check that condition~(\ref{irrotational_2}) is satisfied, and equation~(\ref{QSD_body_general}) may be written as
\begin{eqnarray*}
u_{\bx} &=& \sqrt{\frac{\mbox{det}(\Sigma)}{(2 \pi N)^k}
\left( \frac{b_0 ( y_1^* + y_2^* ) d_0 ( y_1^* + y_2^* ) b_1 ( y_1^* ) d_1 ( y_1^* ) b_2 ( y_2^* ) d_2 ( y_2^* )}{b_0 ( y_1 + y_2 ) d_0 ( y_1 + y_2 ) b_1 ( y_1 ) d_1 ( y_1 ) b_2 ( y_2 ) d_2 ( y_2 )} \right)}
\\ && {} \times
\exp \left( N V ( y_1 , y_2 ) + O (1/N)\right)
\end{eqnarray*}
for $d ( \bx , \delta S ) \ne o(N)$, with $\by = \bx/N$, where $V(y_1, y_2)$ is given by equation~(\ref{V_competition}), and the determinant is given by
\begin{eqnarray}
\mbox{det} ( \Sigma ) &=&
4 
\left( \frac{d_0^\prime (y_1^* + y_2^*)}{d_0 ( y_1^* + y_2^* )}
- \frac{b_0^\prime ( y_1^* + y_2^*)}{b_0 ( y_1^* + y_2^* )} \right)
\left( \frac{d_3^\prime ( y_1^* - y_2^* )}{d_3 ( y_1^* - y_2^* )}
- \frac{b_3^\prime ( y_1^* - y_2^* )}{b_3 ( y_1^* - y_2^* )} \right)
\nonumber \\ && {} +
\left( \frac{d_1^\prime ( y_1^* )}{d_1 ( y_1^* )}
- \frac{b_1^\prime ( y_1^* )}{b_1 ( y_1^* )} \right)
\left( \frac{d_2^\prime ( y_2^* )}{d_2 ( y_2^* )}
- \frac{b_2^\prime ( y_2^* )}{b_2 ( y_2^* )} \right)
\nonumber \\ && {} +
\left( \frac{d_0^\prime (y_1^* + y_2^*)}{d_0 ( y_1^* + y_2^* )}
- \frac{b_0^\prime ( y_1^* + y_2^*)}{b_0 ( y_1^* + y_2^* )}
+ \frac{d_3^\prime ( y_1^* - y_2^* )}{d_3 ( y_1^* - y_2^* )}
- \frac{b_3^\prime ( y_1^* - y_2^* )}{b_3 ( y_1^* - y_2^* )} \right)
\nonumber \\ && {} \times
\left( \frac{d_1^\prime ( y_1^* )}{d_1 ( y_1^* )}
- \frac{b_1^\prime ( y_1^* )}{b_1 ( y_1^* )}
+ \frac{d_2^\prime ( y_2^* )}{d_2 ( y_2^* )}
- \frac{b_2^\prime ( y_2^* )}{b_2 ( y_2^* )} \right) .
\label{det_competition}
\end{eqnarray}
The expression on the right hand side of equation~(\ref{det_competition}) may be slightly simplified by making use of the relationship $b_3 (u) d_3 (u) = a_0$, but we will require the more general form given in equation~(\ref{det_competition}) below.

Taking $a_5 = a_6 = 0$, and correspondingly redefining ${\cal L} = \{ \be_1 , - \be_1 , \be_2 , - \be_2 \}$,  the transition rate functions~(\ref{competition_rates}) are now of the form of the `reversible competition processes' of the corollary to theorem~5 of~\cite{I64}.
From definitions~(\ref{bi_di}) we find
\begin{eqnarray}
\left. \begin{array}{rcl}
b_{11} &=& a_2 b_0^\prime (0) b_1(0) c_2(0) b_3(0), \\
b_{22} &=& a_4 b_0^\prime (0) c_1(0) b_2(0) d_3(0), \\
d_1 &=& a_1 d_0(0) d_1^\prime (0) c_2(0) d_3(0) , \\
d_2 &=& a_3 d_0(0) c_1(0) d_2^\prime (0) b_3(0) .
\end{array} \right\}
\label{bi_di_explicit_rates}
\end{eqnarray}
In order that condition~(\ref{b_d_assumptions}) is satisfied, we require $d_0(0) , b_0^\prime (0) , d_1^\prime (0) , d_2^\prime (0) > 0$.
It is straightforward to check that conditions~(\ref{irrotational_for_BD_process}),~(\ref{irrotational_for_BD_process2}) and~(\ref{linear_Kolmogorov}) are satisfied.
Note that the constraint $b_3(u) d_3(u) = a_0$ is no longer required.
The relationship~(\ref{tau_formula}) now reduces to
\begin{eqnarray*}
\tau &\sim&
\frac{1}{D}
\sqrt{
\frac{2 \pi}{N \mbox{det} ( \Sigma )}
\frac{b_0^\prime (0) d_0(0) b_1(0) d_1^\prime (0) b_2(0) d_2^\prime (0) b_3(0) d_3(0)}{b_0 ( y_1^* y_2^* ) d_0 ( y_1^* + y_2^* ) b_1 ( y_1^* ) d_1 ( y_1^* ) b_2 ( y_2^* ) d_2 ( y_2^* ) b_3 ( y_1^* - y_2^* ) d_3 ( y_1^* - y_2^* )}
}
\\ && {} \times
\exp ( N V (0,0) ) 
\end{eqnarray*}
where $\mbox{det} ( \Sigma )$ is given by equation~(\ref{det_competition}), $D$ is the unique positive solution of equation~(\ref{D_definition}) with $b_{11} , b_{22} , d_1 , d_2$ given by equations~(\ref{bi_di_explicit_rates}), and the function $V ( y_1 , y_2 )$ is given by equation~(\ref{V_competition}).
Condition~(\ref{b_d_assumptions}) ensures that the integrals in equation~(\ref{V_competition}) with $y_1 = y_2 = 0$ all converge.

\section{Derivations}\label{derivations}
We first derive the general results~4 and~5, from which results~1 and~2 follow more or less immediately.
Finally, we derive the more precise result~3.
\subsection{Derivation of result~4}
By~\cite{DP13} theorems~6 and~13 and corollary~7, it follows from our assumption~\ref{QSD_exists} that the process $\bX^{(N)} (t)$ is almost surely absorbed at ${\bf 0}$ within finite time, that the elements of the quasistationary distribution $\bu$ satisfy
\begin{eqnarray}
\sum_{\bl \in {\cal L}} \left( u_{\bx-\bl}  \beta_{\bl} \left( {\bx - \bl \over N} \right) -  u_{\bx} \beta_{\bl} \left( {\bx \over N} \right) \right) &=& - (\tau N)^{-1} u_{\bx} \mbox{ for } \bx \in C , \label{QSD}
\end{eqnarray} 
that the time to extinction for the process initiated according to the quasistationary distribution~$\bu$ is exponentially distributed with mean~$\tau$, and that
\begin{eqnarray}
\tau &=& \left( N  \sum_{\bl \in {\cal L}} u_{- \bl} \beta_{\bl} \left( - \, {\bl \over N} \right) \right)^{-1} . \label{tau_formula}
\end{eqnarray}

Writing $\by = \bx / N$, then, following the methodology described in~\cite{AM17} and references therein, we adopt the  WKB (Wentzel, Kramers, Brillouin) {\em ansatz} that the components $u_{\bx}$ of the quasistationary distribution~$\bu$ may be written in the form~(\ref{ansatz}), where, without loss of generality, we set $V ( \by^* ) = V_0 ( \by^* ) = 0$.
Substituting from~(\ref{ansatz}) into equation~(\ref{QSD}), assuming that $\tau$ is sufficiently large for the right hand side to be neglected, and collecting together terms of leading order in~$N$, we obtain the Hamilton-Jacobi equation
\begin{eqnarray}
\sum_{\bl \in {\cal L}}  \beta_{\bl} ( \by )\left( \exp \left( \bl^T {\partial V \over \partial \by} \right) - 1 \right)
&=& 0 . \label{HJE}
\end{eqnarray}

For later reference, collecting together second order terms gives the transport equation
\begin{eqnarray}
\sum_{\bl \in {\cal L}}  
\exp \left( \bl^T {\partial V \over \partial \by} \right)
\bl^T 
\left(
\left( {\partial V_0 \over \partial \by} 
- {1 \over 2} {\partial^2 V \over \partial \by^2} \bl \right) \beta_{\bl} ( \by )
-  {\partial \beta_{\bl} \over \partial \by}
\right) 
&=& 0 , \label{transport}
\end{eqnarray}
a first-order linear partial differential equation to be solved for $V_0 ( \by )$, once $V ( \by )$ has been found from equation~(\ref{HJE}).

Recalling (assumption~\ref{minus_l}) that if $\bl \in {\cal L}$ then $- \bl \in {\cal L}$, the Hamilton-Jacobi equation~(\ref{HJE}) will be satisfied if we can find a function $V ( \by )$ such that $\btheta = \frac{\partial V}{\partial \by}$ satisfies the system of equations~(\ref{linear_system}).

From assumption~\ref{minus_l}, the system~(\ref{linear_system}) is overdetermined and admits at most one solution $\btheta ( \by )$ for each $\by \in \tilde S^\circ$.
Noting that components of $\bl$ are integer-valued for every $\bl \in {\cal L}$, then condition~(\ref{Kolmogorov0}) ensures that the system~(\ref{linear_system}) is consistent, so that there exists a unique solution $\btheta ( \by )$ for each $\by \in \tilde S^\circ$.


Next, we require that there exist a scalar potential $V ( \by )$ such that $\frac{\partial V}{\partial \by} = \btheta ( \by )$.
But for any $\bl_1 , \bl_2 \in {\cal L}$, we have from equations~(\ref{linear_system}) that
\begin{eqnarray*}
\bl_1^T \frac{\partial \btheta}{\partial \by} \bl_2
&=&
\left( \frac{\partial}{\partial \by} \ln \left( \frac{\beta_{-\bl_1} ( \by )}{\beta_{\bl_1} ( \by )} \right) \right)^T \bl_2 ,
\end{eqnarray*}
so that applying condition~(\ref{irrotational}),
\begin{eqnarray*}
\bl_1^T \frac{\partial \btheta}{\partial \by} \bl_2
&=&
\bl_2^T \frac{\partial \btheta}{\partial \by} \bl_1 .
\end{eqnarray*}
Since ${\cal L}$ spans $\R^k$, it follows that $\frac{\partial \btheta}{\partial \by}$ is a symmetric matrix, 
%
%
%
%
By the Poincar\'e Lemma (\cite{C13}. section~6.17), $\btheta ( \by )$ is therefore the gradient of a scalar potential $V ( \by )$.
The function $V ( \by )$ is given by
\begin{eqnarray}
V ( \by ) &=& \int_{\Gamma} \btheta ( \by^\prime ) \cdot d\by^\prime \mbox{ for } \by \in \tilde S , \label{V_general}
\end{eqnarray}
where $\Gamma$ is any path from $\by^*$ to $\by$ that lies entirely within~$\tilde S$, and the integral is independent of the particular path~$\Gamma$.

From equations~(\ref{tau_formula}) and~(\ref{ansatz}), since ${\cal L}$ is a finite set, we have $( \ln \tau ) / N \to V ( {\bf 0} )$ as $N \to \infty$, and the relationship~(\ref{log_limit_general}) follows, provided the integral converges (noting that the integrand is undefined at the origin).

\subsection{Derivation of result~5}
To evaluate the constant $M_N$ in the WKB expression~(\ref{ansatz}), consider the Taylor series expansion of formula~(\ref{ansatz}) in the vicinity of $\by^*$.
Defining constants
\begin{eqnarray*}
a_{\bl} &=& \bl^T \left. \frac{\partial V}{\partial \by} \right|_{\by = \by^*} \mbox{ for } \bl \in {\cal L} ,
\end{eqnarray*}
then since $\by^*$ is an equilibrium point of the system~(\ref{ODEsystem}) we have 
\begin{eqnarray}
\sum_{\bl \in {\cal L}} \beta_{\bl} ( \by^* ) a_{\bl} = 0 .
\label{equilibrium_V}
\end{eqnarray}
Subtracting equation~(\ref{equilibrium_V}) from equation~(\ref{HJE}) at~$\by^*$, we obtain
\begin{eqnarray}
\sum_{\bl \in {\cal L}} \beta_{\bl} ( \by^* ) \left( \exp \left( a_{\bl} \right) - 1 - a_{\bl} \right) = 0 .
\label{equilibrium_V2}
\end{eqnarray}
Now the function $f(a) = e^a - 1 - a$ is non-negative, and equal to zero only when $a=0$, so that the left hand side of equation~(\ref{equilibrium_V2}) is a sum of non-negative terms, each term of which must therefore be zero, implying that $a_{\bl} = 0$ for all $\bl \in {\cal L}$.
Since ${\cal L}$ spans $\R^k$, it follows that 
\begin{eqnarray}
\left. \frac{\partial V}{\partial \by} \right|_{\by = \by^*} &=& {\bf 0} .
\label{dV_zero}
\end{eqnarray}

Since $\frac{\partial V}{\partial \by} = \btheta ( \by )$, where $\btheta ( \by )$ is the solution of equations~(\ref{linear_system}), we also have
\begin{eqnarray}
\left. \frac{\partial^2 V}{\partial \by^2} \right|_{\by = \by^*}
&=&
\Sigma , \label{Hessian}
\end{eqnarray}
where $\Sigma$ is given by equation~(\ref{Sigma_matrix}).
Recalling the conditions $V ( \by^* ) = V_0 ( \by^* ) = 0$, then for $\left| \bx - N \by^* \right| = O ( \sqrt{N} )$, Taylor series expansion of equation~(\ref{ansatz}) about~$\by^*$ therefore gives
\begin{eqnarray}
u_{\bx} 
&=& M_N \exp \left( - \frac{1}{2N} \left( \bx - N \by^* \right)^T 
\Sigma \left( \bx - N \by^* \right) + o(1) \right) .
\label{Taylor_general}
\end{eqnarray}

Denote by $J ( \by )$ the Jacobian of the system~(\ref{ODEsystem}), so that
\begin{eqnarray*}
J ( \by ) &=& \sum_{\bl \in {\cal L}} \bl \left( \frac{\partial \beta_{\bl}}{\partial \by} \right)^T ,
\end{eqnarray*}
and define the matrix $G$ to be
\begin{eqnarray}
G &=& \sum_{\bl \in {\cal L}} \beta_{\bl} ( \by^* ) \bl \bl^T .
\label{G_definition}
\end{eqnarray}

Then we have
\begin{eqnarray*}
G \Sigma &=& \sum_{\bl \in {\cal L}} \beta_{\bl} ( \by^* ) \bl \bl^T \left. \frac{\partial \btheta}{\partial \by} \right|_{\by = \by^*} \\
&=& \sum_{\bl \in {\cal L}} \beta_{\bl} ( \by^* ) \bl \left( \left. \frac{\partial}{\partial \by} \ln \left( \frac{\beta_{-\bl} ( \by )}{\beta_{\bl} ( \by )} \right) \right|_{\by = \by^*} \right)^T \mbox{ (from equation~(\ref{linear_system}))} \\
&=& \sum_{\bl \in {\cal L}} \beta_{\bl} ( \by^* ) \bl 
\left( 
\frac{1}{\beta_{-\bl} ( \by^* )} \left. \frac{\partial \beta_{-\bl}}{\partial \by} \right|_{\by = \by^*} 
- \frac{1}{\beta_{\bl} ( \by^* )} \left. \frac{\partial \beta_{\bl}}{\partial \by} \right|_{\by = \by^*} \right)^T .
\end{eqnarray*}

Now from equation~(\ref{linear_system}) and equation~(\ref{dV_zero}), we have $\beta_{-\bl} ( \by^* ) = \beta_{\bl} ( \by^* )$ for all $\bl \in {\cal L}$, and so
\begin{eqnarray*}
G \Sigma 
&=& \sum_{\bl \in {\cal L}} \bl 
\left( 
\left. \frac{\partial \beta_{-\bl}}{\partial \by} \right|_{\by = \by^*} 
-  \left. \frac{\partial \beta_{\bl}}{\partial \by} \right|_{\by = \by^*} \right)^T .
\end{eqnarray*}
But since $\bl \in {\cal L} \Leftrightarrow - \bl \in {\cal L}$, this can be written as
\begin{eqnarray}
G \Sigma 
&=& - \sum_{\bl \in {\cal L}} - \bl 
\left( \left. \frac{\partial \beta_{-\bl}}{\partial \by} \right|_{\by = \by^*} \right)^T
-  \sum_{\bl \in {\cal L}} \bl \left( \left. \frac{\partial \beta_{\bl}}{\partial \by} \right|_{\by = \by^*} \right)^T 
\nonumber \\
&=& - 2 J ( \by^* ) .
\label{Sigma_equation}
\end{eqnarray}

It is clear from the definition~(\ref{G_definition}) that $G$ is symmetric and positive-definite, and from the equation~(\ref{Hessian}) that $\Sigma$ is symmetric.
From assumption~\ref{equilm_assumptions} we know that all eigenvalues of $J ( \by^* )$ have strictly negative real part, so it follows from equation~(\ref{Sigma_equation}) that $\Sigma$ is invertible, and that $\Sigma^{-1}$ satisfies the Lyapunov equation $J ( \by^* ) \Sigma^{-1} + \Sigma^{-1} J ( \by^* )^T + G = 0$.
From theorems~13.21 and~13.24 of~\cite{L05}, this Lyapunov equation has a unique solution $\Sigma^{-1}$, and $\Sigma^{-1}$ is positive definite.
Hence equation~(\ref{Taylor_general}) represents a multivariate Gaussian distribution with variance matrix $N \Sigma^{-1}$, normalisation of which implies that
\begin{eqnarray}
M_N &=& \sqrt{\frac{\mbox{det} ( \Sigma )}{(2 \pi N )^k}} .
\label{M_N}
\end{eqnarray}

Now consider the function $V_0 ( \by )$.
Since $\btheta ( \by ) = \frac{\partial V}{\partial \by}$ satisfies equations~(\ref{linear_system}), the transport equation~(\ref{transport}) simplifies to
\begin{eqnarray*}
\sum_{\bl \in {\cal L}} 
\beta_{-\bl} ( \by ) \bl^T 
\left( \frac{\partial V_0}{\partial \by} 
- \frac{1}{2} \frac{\partial}{\partial \by} \ln \left( \beta_{-\bl} ( \by ) \beta_{\bl} ( \by ) \right) \right) 
&=& 0 . 
\end{eqnarray*}

This will be satisfied if we can find a function $V_0 ( \by )$ such that $\btheta^{0} = \frac{\partial V_0}{\partial \by}$ satisfies the system of equations~(\ref{linear_system2}).
Arguing as in the derivation of result~4, condition~(\ref{Kolmogorov1}) ensures that the linear system~(\ref{linear_system2}) admits a unique solution $\btheta^{0} ( \by )$, and condition~(\ref{irrotational_2}) ensures, via the Poincar\'e Lemma, that $\btheta^{0} ( \by )$ is the gradient of a scalar potential $V_0 ( \by )$, with
\begin{eqnarray}
V_0 ( \by ) &=& \int_{\Gamma^\prime} \btheta^{0} ( \by^\prime ) \cdot d \by^\prime \mbox{ for } \by \in \tilde S^\circ ,
\label{V0_integral}
\end{eqnarray}
where $\Gamma^\prime$ is any path from $\by^*$ to $\by$ that lies entirely within $\tilde S$, and the integral is independent of the particular path~$\Gamma^\prime$.

Since ${\cal L}$ spans $\R^k$, any $\by \in \tilde S$ may be written in the form $\by = \by^* + \sum_{i=1}^n a_i \bl_i$ for some (not necessarily unique) $a_1 , a_2 , \ldots , a_n \in \R$ and $\bl_1 , \bl_2 , \ldots , \bl_n \in {\cal L}$ such that $\by^* + \sum_{i=1}^j a_i \bl_i \in \tilde S$ for $j=1,2,\ldots,n$.
From equations~(\ref{linear_system2}), we have
\begin{eqnarray*}
\bl^T \frac{\partial}{\partial \by} \left( V_0 - \frac{1}{2} \ln \left( \beta_{-\bl} ( \by ) \beta_{\bl } ( \by ) \right) \right) &=& 0
\mbox{ for all } \bl \in {\cal L} .
\end{eqnarray*}
Starting from $\by^*$ and integrating along each direction $\bl_1 , \bl_2 , \ldots , \bl_n$ in turn, we thus obtain
\begin{eqnarray}
V_0 ( \by ) &=& 
\frac{1}{2} \sum_{j=1}^n \ln \left( 
\frac
{\beta_{-\bl_j} \left( \by^* + \sum_{i=1}^j a_i \bl_i \right) 
\beta_{\bl_j} \left( \by^* + \sum_{i=1}^j a_i \bl_i \right)}
{\beta_{-\bl_j} \left( \by^* + \sum_{i=1}^{j-1} a_i \bl_i \right) 
\beta_{\bl_j} \left( \by^* + \sum_{i=1}^{j-1} a_i \bl_i \right)}
\right) .
\label{V0_general}
\end{eqnarray}
Since the integral in equation~(\ref{V0_integral}) is independent of the particular path chosen, the expression~(\ref{V0_general}) is correspondingly independent of the particular representation $\by = \by^* + \sum_{i=1}^n a_i \bl_i$.

Notice that the expression~(\ref{V0_general}) diverges for $\by \in \partial \tilde S$.
For $d ( \by , \partial \tilde S ) \ne o(1)$, substituting for $V ( \by )$, $V_0 ( \by )$ and $M_N$ from equations~(\ref{V_general}),~(\ref{V0_general}) and~(\ref{M_N}) into the WKB formula~(\ref{ansatz}), the result follows.

\subsection{Derivation of result~1}
In the case ${\cal L} = \{ \be_i , - \be_i : i=1,2,\ldots,k \}$, condition~(\ref{Kolmogorov0}) is automatically satisfied, and condition~(\ref{irrotational}) reduces to condition~(\ref{irrotational_for_BD_process}).
The system of equations~(\ref{linear_system}) reduces to $\theta_i ( \by ) = h_i ( \by )$ for $i=1,2,\ldots,k$, where $h_i ( \by )$ is given by equations~(\ref{h_components_BD}), and so the function $V ( \by )$ in the derivation of result~4 is given by
\begin{eqnarray}
V ( \by ) &=& \int_{\Gamma} \bh ( \by^\prime ) \cdot d\by^\prime \mbox{ for } \by \in \tilde S \label{V_for_BD_process}
\end{eqnarray}
where $\Gamma$ is any path from $\by^*$ to $\by$ that lies entirely within~$\tilde S$, and the integral is independent of the particular path~$\Gamma$.

The integrand in equation~(\ref{V_for_BD_process}) is undefined at the origin.
By assumptions~\ref{absorbing} and~\ref{communication}, recalling the definitions~(\ref{bi_di}), we have that $b_{ij} \ge 0$, $d_i \ge 0$ for all $i,j$, and that $\left. \frac{\partial \beta_{-\be_i}}{\partial y_j} \right|_{\by = {\bf 0}} = 0$ for $i \ne j$.
We now make the further supposition, condition~(\ref{b_d_assumptions}), that $d_i > 0$ for all~$i$, and that for each~$i$, there exists some~$j$ such that $b_{ij} > 0$.
Consider approach to the origin along the path $\by = \hat y \bxi$ for some fixed $\bxi = \left( \xi_1 , \xi_2 , \ldots , \xi_k \right)$ with $\xi_i > 0$ for all~$i$ and $\sum_{i=1}^k \xi_i = 1$.
Applying l'H\^opital's rule along this path, we have
\begin{eqnarray*}
\lim_{\hat y \to 0} \frac{\beta_{-\be_i} ( \hat y \bxi )}{\beta_{\be_i} ( \hat y \bxi )} &=&
\frac{\xi_i d_i}{\sum_{j=1}^k \xi_j b_{ij}} ,
\end{eqnarray*}
a finite and non-zero limit. 
Hence along this path, for each~$i$, $h_i ( \by )$ defined by~(\ref{h_components_BD}) converges to a finite limit, and hence for $\by = {\bf 0}$, the integral in equation~(\ref{V_for_BD_process}) converges.

Result~1 now follows from result~4.

\subsection{Derivation of result~2}
In the case ${\cal L} = \{ \be_i , - \be_i : i=1,2,\ldots,k \}$, condition~(\ref{Kolmogorov1}) is automatically satisfied, and condition~(\ref{irrotational_2}) reduces to condition~(\ref{irrotational_for_BD_process2}).
The system of equations~(\ref{linear_system2}) reduces to $\theta_i^{0} ( \by ) = h_i^0 ( \by )$ for $i=1,2,\ldots,k$, where $h_i^0 ( \by )$ is given by equations~(\ref{h0_BD}), and result~2 thus follows from result~5.

\subsection{Derivation of result~3}
To approximate the mean extinction time~$\tau$ using formula~(\ref{tau_formula}), we need to be able to approximate the quasistationary probabilities $u_{\be_i}$ for $i=1,2,\ldots,k$.
Since $V_0 ( \by )$ diverges for $\by \in \partial \tilde S$, we cannot directly make use of the WKB approximation~(\ref{QSD_body}).
Instead, we will now seek an approximation for $u_{\bx}$ valid for $| \bx | = O(1)$.
In order to normalise this approximation, we will make use of the WKB~approximation~(\ref{QSD_body}).

With the convention that $u_{\bx} = 0$ for $\bx \notin C$, the exact equation~(\ref{QSD}) may be written as
\begin{eqnarray}
\sum_{i=1}^k \left( 
u_{\bx - \be_i} \beta_{\be_i} \left( \frac{\bx - \be_i}{N} \right) 
+ u_{\bx + \be_i} \beta_{-\be_i} \left( \frac{\bx + \be_i}{N} \right)
- u_{\bx} \left( \beta_{\be_i} \left( \frac{\bx}{N} \right) + \beta_{-\be_i} \left( \frac{\bx}{N} \right) \right)
\right) \nonumber
\\
&&\hspace*{-50mm} = \;\; - ( \tau N )^{-1} u_{\bx} \mbox{ for } \bx \in C . \label{BD_balance}
\end{eqnarray}

Assuming, as before, that $\tau$ is sufficiently large for the right-hand side of equation~(\ref{BD_balance}) to be neglected, taking the linear approximation to the left-hand side, recalling the definitions~(\ref{bi_di}) of the constants $b_{ij} , d_i$, and writing $\bm{b}^{(i)} = \left( b_{i1} , b_{i2} , \ldots , b_{ik} \right)$ for $i=1,2,\ldots,k$, we obtain the asymptotic balance equation
\begin{eqnarray}
\sum_{i=1}^k \left( 
\tilde u_{\bx - \be_i} ( \bx - \be_i )^T 
\bm{b}^{(i)}
+ \tilde u_{\bx + \be_i} ( x_i + 1 ) 
d_i
\right)
- \tilde u_{\bx} \sum_{i=1}^k \left( \bx^T 
\bm{b}^{(i)}
+ x_i d_i \right)
&=& 0 , \label{linear_balance}
\end{eqnarray}
where $\{ \tilde u_{\bx} : \bx \in C \}$ denotes the solution of the linearised balance equation.

Recalling our assumptions that  $b_{ii} , d_i > 0$ for $i=1,2,\ldots,k$, then under the condition~(\ref{linear_Kolmogorov}), which ensures that the linear approximating process satisfies the Kolmogorov condition~(\ref{Kolmogorov}), the asymptotic balance equation~(\ref{linear_balance}) has solution
\begin{eqnarray}
\tilde u_{\bx}
&=& \Lambda \,
\frac{1}{\sum_{i=1}^k x_i} \frac{\left( \sum_{i=1}^k x_i \right) !}{\prod_{i=1}^k x_i !} \prod_{i=1}^k b_{ii}^{x_i}
\left(
\prod_{i=1}^k \left( d_i \right)^{-x_i}
-
\prod_{i=1}^k \left( D + d_i \right)^{-x_i}
\right) 
\label{u_tilde}
\end{eqnarray}
where $\Lambda$ is a normalising constant to be found, and $D$ satisfies equation~(\ref{D_definition}).
The first component of the solution~(\ref{u_tilde}) is found by solving the detailed balance equations~(\ref{detailed_balance}) corresponding to equation~(\ref{linear_balance}), and the other component by analogy with the solution for the case $k = 1$ given in~\cite{AM10}, see also~\cite{C18b}.

Local instability of the equilibrium point at $\by = {\bf 0}$ for the system~(\ref{ODEsystem}) (assumption~\ref{equilm_assumptions}) implies that at least one eigenvalue of the Jacobian of system~(\ref{ODEsystem}) at the origin, $J ( {\bf 0} )$, has positive real part.
Now $J ( { \bf 0} )$ has components
\begin{eqnarray*}
J_{ij} ( {\bf 0} ) &=& b_{ii} - d_i \delta_{ij} \mbox{ for } i,j = 1,2,\ldots,k ,
\end{eqnarray*}
where $\delta_{ij}$ is the Kronecker delta.
It follows from theorem~A.1 of~\cite{DHR10} that $\sum_{i=1}^k \left( b_{ii} / d_i \right) > 1$, and hence that equation~(\ref{D_definition}) has a unique positive solution~$D$.

To evaluate the normalising constant~$\Lambda$ in formula~(\ref{u_tilde}), we match the approximation~(\ref{u_tilde}) with our WKB approximation~(\ref{QSD_body}) in their common range of validity.
First, consider the approximation~(\ref{u_tilde}) for large $| \bx |$.
Applying Stirling's formula to the factorial terms, and recalling that $D>0$, we obtain
\begin{eqnarray*}
\tilde u_{\bx} &\sim&
\Lambda \,
\left( \sum_{i=1}^k x_i \right)^{\sum_{i=1}^k x_i}
\sqrt{\frac{1}{(2 \pi)^{k-1} \left( \sum_{i=1}^k x_i \right) \prod_{i=1}^k x_i}}
\,
\prod_{i=1}^k \left( \frac{b_{ii}}{x_i d_i} \right)^{x_i} .
\end{eqnarray*}

Along the specific trajectory $\bx = \hat x \bxi$, where $\bxi = \left( \xi_1 , \xi_2 , \ldots , \xi_k \right)$ is fixed, with $\xi_1 , \xi_2 , \ldots , \xi_k > 0$ and $\sum_{i=1}^k \xi_i = 1$, then as $\hat x \to \infty$ we have
\begin{eqnarray}
\tilde u_{\bx} &\sim&
\Lambda \,
\sqrt{\frac{1}{(2 \pi)^{k-1} \hat x^{k+1}  \prod_{i=1}^k \xi_i}}
\, \prod_{i=1}^k \left( \frac{b_{ii}}{\xi_i d_i} \right)^{\hat x \xi_i} .
\label{linear_asymptote}
\end{eqnarray}

Returning to the WKB approximation~(\ref{QSD_body}), consider the trajectory $\by = \hat y \bxi$ as $\hat y \to 0$, where $\by = \bx / N$.
Under the assumption~(\ref{linear_Kolmogorov}), applying l'H\^opital's rule along the trajectory $\by = \hat y \bxi$ to the derivatives $\left. \partial V \right/ \partial y_i = h_i ( \by )$ given by~(\ref{h_components_BD}), we have
\begin{eqnarray*}
\frac{\partial V}{\partial y_i} &=& \ln \left( \frac{\beta_{-\be_i} ( \hat y \bxi )}{\beta_{\be_i} ( \hat y \bxi )} \right) 
\to 
\ln \left( 
\frac{\xi_i d_i}{b_{ii}} 
\right)
\mbox{ as } \hat y \to 0 ,
\end{eqnarray*}
so that Taylor series expansion of $V ( \by )$ in the range $| \bx | = O ( \sqrt{N} )$ gives
\begin{eqnarray}
\exp ( - N V ( \by ) ) 
&\sim&
\exp \left( - N V ( {\bf 0} )
- \hat x \sum_{i=1}^k \xi_i \ln \left( 
\frac{\xi_i d_i}{b_{ii}} 
\right)
+ O (1/N) \right) \nonumber \\
&=& 
\prod_{i=1}^k \left( 
\frac{b_{ii}}{\xi_i d_i}
\right)^{\hat x \xi_i}
\exp \left( - N V ( {\bf 0} )
+ O (1/N) \right) . \label{Taylor_V}
\end{eqnarray}

Under the assumption~(\ref{linear_Kolmogorov}), we also have that for $i=1,2,\ldots,k$,
\begin{eqnarray}
\beta_{-\be_i} ( y_1 , \ldots , y_i , y_{i+1}^* , \ldots , y_k^* )
&=&
\hat y \xi_i \left. 
\frac{\partial \beta_{-\be_i}}{\partial y_i} 
\right|_{\by = (0,\ldots,0,y_{i+1}^*, \ldots, y_k^* )}
+ O ( \hat y^2 ) , \label{Taylor_death} 
\end{eqnarray}
and
\begin{eqnarray}
\beta_{\be_k} ( {\by} )
&=& 
b_{kk} \hat y + O ( \hat y^2 ) . \label{Taylor_birth}
\end{eqnarray}

Substituting from equations~(\ref{Taylor_V},~\ref{Taylor_death},~\ref{Taylor_birth}) into the WKB approximation formula~(\ref{QSD_body}) yields, along the trajectory $\bx = \hat x \bxi$ in the range $| \bx | = O ( \sqrt{N} )$, 
\begin{eqnarray}
\tilde u_{\bx} 
&\sim&
\sqrt{\frac{N \mbox{det} ( \Sigma )}{(2 \pi)^k \hat x^{k+1}}
{\prod_{i=1}^k
\beta_{\be_i} \left( 0 , \ldots , 0 , y_i^* , \ldots , y_k^* \right) \beta_{-\be_i} \left( 0 , \ldots , 0 , y_i^* , \ldots , y_k^* \right)
\over
b_{kk}
\prod_{i=1}^k \left( \xi_i \left. \frac{\partial \beta_{-\be_i}}{\partial y_i} \right|_{\by = (0 , \ldots , 0 , y^*_{i+1} , \ldots , y^*_k )} \right)
\prod_{i=1}^{k-1} \beta_{\be_i} \left( 0 , \ldots , 0 , y_{i+1}^* , \ldots , y_k^* \right) 
}}
\nonumber \\ && {} \times
\prod_{i=1}^k \left( 
\frac{b_{ii}}{\xi_i d_i}
\right)^{\hat x \xi_i}
\exp \left( - N V ( {\bf 0} ) 
\right) . \label{WKB_asymptote}
\end{eqnarray}

Matching the expressions~(\ref{linear_asymptote}) and~(\ref{WKB_asymptote}), we see that terms in $\hat x$ and $( \xi_1 , \xi_2 , \ldots , \xi_k )$ do indeed match, and we obtain the normalising constant~$\Lambda$ as
\begin{eqnarray}
\Lambda
&=&
\sqrt{\frac{N \mbox{det} ( \Sigma )}{2 \pi}
{\prod_{i=1}^k
\beta_{\be_i} \left( 0 , \ldots , 0 , y_i^* , \ldots , y_k^* \right) \beta_{-\be_i} \left( 0 , \ldots , 0 , y_i^* , \ldots , y_k^* \right)
\over
b_{kk}
\prod_{i=1}^k \left. \frac{\partial \beta_{-\be_i}}{\partial y_i} \right|_{\by = (0 , \ldots , 0 , y^*_{i+1} , \ldots , y^*_k )}
\prod_{i=1}^{k-1} \beta_{\be_i} \left( 0 , \ldots , 0 , y_{i+1}^* , \ldots , y_k^* \right) 
}}
\nonumber \\ 
&& {} \times
\exp \left( - N V( {\bf 0} ) \right) .
\label{A_formula}
\end{eqnarray}

Now substituting for $\tilde u_{\be_i}$ from equation~(\ref{u_tilde}) into equation~(\ref{tau_formula}), noting that $\beta_{-\be_i} \left( \frac{\be_i}{N} \right) = d_i + O(1/N)$, and recalling the relationship~(\ref{D_definition}), we obtain
\begin{eqnarray}
\tau &=& \left( N \sum_{i=1}^k u_{\be_i} \beta_{-\be_i} \left( \frac{\be_i}{N} \right) \right)^{-1} \nonumber
\\
&\sim&
\left( \Lambda \sum_{i=1}^k b_{ii} \left( 1 - \frac{d_i} {D + d_i} \right) 
\right)^{-1} \nonumber
\\
&=& \left( \Lambda D \right)^{-1} . \label{tau_AD}
\end{eqnarray}
From equations~(\ref{tau_AD}),~(\ref{A_formula}) and~ (\ref{V_for_BD_process}), the result~(\ref{tau_BD}) follows.

{\bf Acknowledgements.}
It is a pleasure to thank Lyonell Boulton and Heiko Gimperlein for helpful discussions.

\bibliographystyle{plain}
\bibliography{arXiv_version}

\end{document}